\documentclass[fleqn]{mat01}
\usepackage{times,mathtimy,amssymb}
\begin{document}

\setcounter{page}{243}
\firstpage{243}

\def\theore{\trivlist\item[\hskip\labelsep{{\bf Theorem.}}]}

\font\xx=msam5 at 10pt
\def\ab{\mbox{\xx{\char'03}}}

\newtheorem{corolla}{\rm COROLLARY}

\newtheorem{exampl}{Example}

\title{Necessary and sufficient conditions for inclusion relations\\ for
absolute summability}

\markboth{B~E~Rhoades and Ekrem Sava\c{s}}{Inclusion relations}

\author{B~E~RHOADES and EKREM SAVA\c{S}$^{*}$}

\address{Department of Mathematics, Indiana University,
Bloomington, IN 47405-7106, USA\\
\noindent $^{*}$Department of Mathematics, Y\"{u}z\"{u}nc\"{u} Yil University,
Van, Turkey\\
\noindent E-mail: rhoades@indiana.edu; ekremsavas@yahoo.com}

\volume{113}

\mon{August}

\parts{3}

\Date{MS received 17 January 2001; revised 19 August 2002}

\begin{abstract}
We obtain a set of necessary and sufficient conditions for $| 
\overline{N}, p_{n} |_{k} $ to imply $|\overline{N}, q_{n} |_{s}$ for $1 < k \leq
s < \infty$. Using this result we establish several inclusion theorems
as well as conditions for the equivalence of $| \overline {N}, p_{n}
|_{k}$ and $|\overline{N},q_{n} |_{s}$.
\end{abstract}

\keyword{Absolute summability; weighted mean matrix; Ces\'{a}ro matrix.}

\maketitle\vspace{1.5pc}

\noindent In 1994, Sar{\i}g\"{o}l \cite{Sa94} obtained necessary and sufficient
conditions for $| {\overline{N},p_{n}} |_{k}$ to imply 
$| {\overline {N},q_{n}} |_{s}$ for $1 < k \le s < 
\infty$, using the definition that a series $\sum{a_{k}}$ is summable in
$| {\overline {N},p_{n}} |_{k}$ iff
\begin{equation}
\sum_{n = 1}^{\infty}  \left( \frac{P_{n}} {p_{n}}
\right) ^{k - 1}\left| \Delta T_{n}  \right|^{k} < \infty,
\end{equation}
where $T_{n}$ is the $n$th term of the ${\left| {\overline
{N},p_{n}}  \right|}$ transform of the sequence of partial sums of the
series $\sum a_{n}$.

As pointed out by the first author in \cite{Rh99} the correct
condition is
\begin{equation}
\sum\limits_{n = 1}^{\infty}  n ^{k - 1}\left| {\Delta T_{n}}
\right|^{k} < \infty.
\end{equation}

In this paper we obtain appropriate necessary and sufficient conditions
for $| \overline{N},p_{n}|_{k}$ summability to imply that of $|
\overline{N}, q_{n} |_{s}$ for $1 < k \leq s < \infty $. As in
\cite{Sa94} we make use of a result of Bennett \cite{Be87}, who has obtained necessary
and sufficient conditions for a factorable matrix to map $\ell ^{k} \to
\ell ^{s}$. A factorable matrix $A$ is one in which each entry $a_{nk} =
b_{n} c_{k}$. Weighted mean matrices are factorable.

It will not be possible to extended our result by replacing $(\overline{N},
q_{n})$ by a triangular matrix $A$, since necessary and sufficient
conditions are not known for an arbitrary triangular matrix $B$ to map:
$\ell ^{k} \to \ell ^{s}$.

However, if $k=1$, then the necessary and sufficient conditions can be
obtained. Such a result is the special case, by setting each
$\lambda_{n} = 1$ of Theorem~2.1 of Rhoades and\break Sava\c{s}
\cite{RS02}.\pagebreak

Our main result is the following:

\begin{theore}
{\it Let ${\left\{ {p_{n}} \right\}}$ and ${\left\{ {q_{n}} \right\}}$ be 
positive sequences{\rm ,} $1<k\le s<\infty$. Then $| 
{\overline{N}, p_{n}} |_{k} \Rightarrow | {\overline{N},
q_{n}} |_{s}$ iff\\[.5pc]

\noindent {\rm (i)}\vspace{-2.17pc}
\begin{equation}
n^{(1/k - 1/s)} \frac{q_{n} P_{n}}{p_{n} Q_{n}} = O(1)
\end{equation}
and\\[.5pc]

\noindent {\rm (ii)}\vspace{-2.65pc}
\begin{equation}
\Biggl( \sum_{n = m}^{\infty} \biggl( n^{1 - 1/s} \frac{q_{n}}{Q_{n}
Q_{n - 1}} \biggl) ^{s} \Biggr)^{1/s}
\left( \sum_{v = 1}^{m} \left| Q_{v} - \frac{q_{v} 
P_{v}}{p_{v}} \right| ^{k^{\ast}}\left( \frac{1}{v} \right)
\right)^{1/k^{\ast}} = O(1),
\end{equation}
\noindent where $k^{*}$ denotes the conjugate index of $k$ i.e., $1/k +
1/k^{*} = 1$.}
\end{theore}

\begin{proof}
Let $(x_{n})$ and $(y_{n})$ denote the $n$th terms of the $|
{\overline{N}, p_{n}} |$ and ${\left| {\overline{N}, q_{n}}
\right|}$, transforms respectively of $s_{n} = {\sum_{i = 0}^{n}
{a_{i}}}$. Then
\begin{equation}
X_{n} := x_{n} - x_{n - 1} = \frac{p_{n}}{P_{n} P_{n - 1}}\sum_{v =
1}^{n-1} P_{v - 1} a_{v}
\end{equation}
and
\begin{equation}
Y_{n} := y_{n} - y_{n - 1} = \frac{q_{n}}{Q_{n} Q_{n - 1}}\sum_{v =
1}^{n} Q_{v - 1} a_{v}.
\end{equation}

Solving~(5) for $a_{n}$ gives
\begin{align*}
&{\frac{{P_{n} P_{n - 1} X_{n}} }{{p_{n}} }} = {\sum\limits_{v = 1}^{n}
{P_{v - 1} a_{v}}},\\[.2pc]
&{\frac{{P_{n - 1} P_{n - 2} X_{n - 1}} }{{p_{n - 1}} }} =
{\sum\limits_{v = 1}^{n} {P_{v - 1} a_{v}}}.
\end{align*}

Thus
\begin{equation*}
\frac{P_{n} P_{n - 1} X_{n}}{p_{n}} - \frac{P_{n - 1} P_{n - 2} X_{n -
1}}{p_{n - 1}} = P_{n - 1} a_{n},
\end{equation*}
or
\begin{equation}
a_{n} = \frac{P_{n} X_{n}}{p_{n}} - \frac{P_{n - 2} X_{n - 1} }{p_{n -
1}}.
\end{equation}

Substituting (7) into (6) we have
\begin{align*}
Y_{n} &= {\frac{{q_{n}} }{{Q_{n} Q_{n - 1}} }}{\sum\limits_{v = 1}^{n}
{Q_{v - 1}} } {\left[ {{\frac{{P_{v} X_{v}} }{{p_{v}} }} - {\frac{{P_{v
- 2} X_{v - 1}} }{{p_{v - 1}} }}} \right]}\\[.2pc]
&= {\frac{{q_{n}} }{{Q_{n} Q_{n - 1}} }}{\left[ {{\sum\limits_{v =
1}^{n} {{\frac{{Q_{v - 1} P_{v} X_{v}} }{{p_{v}} }}}} - {\sum\limits_{v
= 1}^{n} {{\frac{{Q_{v - 1} P_{v - 2} X_{v - 1}} }{{p_{v - 1}} }}}} }
\right]}
\end{align*}
\begin{align*}
&= {\frac{{q_{n}} }{{Q_{n} Q_{n - 1}} }}{\left[ {{\sum\limits_{v =
1}^{n} {{\frac{{Q_{v - 1} P_{v} X_{v}} }{{p_{v}} }}}} - {\sum\limits_{i
= 0}^{n - 1} {{\frac{{Q_{i} P_{i - 1} X_{i}} }{{p_{i}} }}}} } \right]}\\[.2pc]
& = {\frac{{q_{n}} }{{Q_{n} Q_{n - 1}} }}{\left[ {{\frac{{Q_{n - 1}
P_{n} X_{n}} }{{p_{n}} }} + {\sum\limits_{v = 1}^{n - 1} {\left( {P_{v}
Q_{v - 1} - Q_{v} P_{v - 1}} \right){\frac{{X_{v}} }{{p_{v}} }}}} }
\right]}.
\end{align*}

But
\begin{align*}
P_{v} Q_{v - 1} - Q_{v} P_{v - 1} &= P_{v} (Q_{v} - q_{v} ) - Q_{v}
(P_{v} - p_{v} ) \\[.2pc]
&= - q_{v} P_{v} + p_{v} Q_{v}.
\end{align*}

Therefore
\begin{equation*}
Y_{n} = {\frac{{q_{n} P_{n} X_{n}} }{{p_{n} Q_{n}} }} + {\frac{{q_{n}
}}{{Q_{n} Q_{n - 1}} }}{\sum\limits_{v = 1}^{n - 1} {\left( {Q_{v} -
{\frac{{q_{v} P_{v}} }{{p_{v}} }}} \right)}} X_{v}.
\end{equation*}

Define
\begin{equation*}
Y_{n} ^{\ast} = n^{1 - 1/s }Y_{n}, \quad X_{n} ^{\ast} = n^{1 -
1/k}X_{n}.
\end{equation*}

Then
\begin{align*}
Y_{n}^{\ast} &= n^{1 - 1/ s}\Biggl[ \frac{q_{n} P_{n}}{p_{n}
Q_{n}}\left( \frac{X_{n}^{\ast }}{n^{1 - 1/k}} \right) \\[.2pc]
&\quad\ + \frac{q_{n}}{Q_{n} Q_{n - 1}} \sum_{v = 1}^{n - 1} \left(
Q_{v} - \frac{q_{v} P_{v}}{p_{v}} \right) \left( \frac{X_{v}
^{\ast}}{v^{1 - 1/k}} \right) \Biggr].
\end{align*}

Thus, $Y_{n} ^{\ast}  = {\sum_{v = 1}^{n} {a_{nv} X_{v} ^{\ast}}}$, where
\begin{equation}
a_{nv} = \begin{cases}\displaystyle
\frac{n^{1 - 1/s}q_{n}}{Q_{n}Q_{n - 1}}\frac{(Q_{\nu} - \frac{q_{\nu}
P_{\nu}}{p_{\nu}})}{\nu^{1 - 1/k}}, &1 \leq \nu < n\\[1pc]
\displaystyle \frac{n^{1/k-1/s}q_{n}P_{n}}{p_{n}Q_{n}}, &\nu = n\\[.6pc]
0, &\nu > n
\end{cases}.
\end{equation}

Then $| {\overline {N},p_{n}} |_{k} \Rightarrow |
{\overline {N},q_{n}} |_{s} $ is equivalent to
\begin{equation*}
\sum \nolimits | X_n^{\ast} |^k < \infty \quad \Rightarrow \quad
\sum\nolimits | Y_{n}^{\ast} |^{s} < \infty; \quad \textnormal{i.e.,}
\quad A:\ell^{k} \to \ell^{s},
\end{equation*}
where $A$ is the matrix whose entries are defined by (8). We may write $A
= B + C$, where $b_{nv} = a_{nv}$ for $1 \le v < n; b_{nv} = 0$,
otherwise, and $C$ is the diagonal matrix with $c_{nn} = a_{nn}$.
Omitting the first row of $B$, which contains all zeros, what remains is
a factorable matrix.

From Theorem 2(ii) of \cite{Be87} a factorable
matrix with nonzero entries $b_{n} c_{v} $, is a bounded operator
from $\ell ^{k}$ to $ \ell ^{s}$ iff
\begin{equation}
\left( \sum\limits_{n = m}^{\infty} b_{n}^{s}
\right)^{1/s} \left( \sum\limits_{v = 1}^{m} c_{v}^{k^{\ast}}
\right)^{1/k^{*}} = O(1),
\end{equation}
where $k^{\ast}$ is conjugate index to $k$.

\pagebreak

Applying (9) to $B$, and using  (8),  we have that $B:\ell ^{k} \to
\ell ^{s}$ iff
\begin{equation}
\Biggl( \sum_{n = m + 1}^{\infty} \Biggl( \frac{n^{1 - 1/s}
q_{n}}{Q_{n} Q_{n - 1}} \Biggr)^{s} \Biggr)^{1/s}
\left( \sum_{v = 1}^{m}\left| Q_{v} - \frac{q_{v} P_{v}}{p_{v}}
\right|^{k^{\ast}}\left(\frac{1}{v} \right) \right)^{1/k^{\ast}} = O(1).
\end{equation}

Since $k \le s, C:\ell ^{k} \to \ell ^{s}$, i.e.,
\begin{equation}
\left( {{\sum\limits_{n = 1}^{\infty} {{\left| {c{}_{nn}s_{n}}
\right|}^{s}}}} \right)^{1/s} < \infty
\end{equation}
for every ${\left\{ {s_{n}} \right\}} \in \ell ^{k}$. But (11) implies
that ${\left\{ {c_{nn}} \right\}} \in \ell ^{s^{\ast} }$, where
$s^{\ast} $ is the conjugate of $s$. In particular ${\left\{ {c_{nn}}
\right\}}$ is bounded. Conversely if ${\left\{ {c_{nn}} \right\}}$ is
bounded, since $k \le s$, $C:\ell ^{k} \to \ell ^{s}$.

Combining these facts, $A:\ell ^{k} \to \ell ^{s}$ iff (3) and (4) are
satisfied.

This completes the proof of the theorem.\hfill $\ab$
\end{proof}

\begin{corolla}$\left.\right.$\vspace{.5pc}

\noindent Let $\{ q_{n}\} $ be a positive sequence{\rm ,} $1 < k \leq s <
\infty$. Then $| C,1 |_{k} \Rightarrow | \overline {N},q_{n} |_{s} $
iff\\[.5pc]

\noindent {\rm (i)}\vspace{-2.17pc}
\begin{equation*}
\Biggl( \frac{n^{1 + 1/k - 1/s}q_{n}}{Q_{n}} \Biggr) = O(1)
\end{equation*}
and\\[.5pc]

\noindent {\rm (ii)}\vspace{-2.25pc}
\begin{equation*}
\Biggl( \sum_{n = m + 1}^{\infty} \Biggl( \frac{n^{1 - 1/s }q_{n}
}{Q_{n} Q_{n - 1}} \Biggr)^{s} \Biggr)^{1/s}
\biggl( \sum_{v = 1}^{m} | Q_{v} - (v + 1)q_{v} |^{k^{\ast}}
\biggl( \frac{1}{v} \biggr) \biggr)^{1/k^{\ast}} = O(1).
\end{equation*}
\end{corolla}

\begin{proof}
Set $p_{n} \equiv 1$ in Theorem~1.\hfill $\ab$
\end{proof}

\begin{corolla}$\left.\right.$\vspace{.5pc}

\noindent Let ${\left\{ {p_{n}} \right\}}$ be a positive sequence{\rm ,} $1 < k
\le s < \infty$. Then $| \overline{N},p_{n} |_{k}
\Rightarrow {\left| {C,1} \right|}_{s}$ iff\\[.5pc]

\noindent {\rm (i)}\vspace{-2.3pc}
\begin{equation*}
\frac{n^{(1/k)-(1/s)-1} P_{n}}{p_{n}} = O(1)
\end{equation*}
and\\[.5pc]

\noindent {\rm (ii)}\vspace{-2.55pc}
\begin{equation*}
\left( \sum\limits_{v = 1}^{m} \left| v + 1 - \frac{P_{v}}{p_{v}}
\right|^{k^{\ast}} \left( \frac{1}{v} \right) \right)^{1/k^{\ast}} = O(m).
\end{equation*}
\end{corolla}

\begin{proof} 
In Theorem~1, set $q_{n} \equiv 1$, to obtain condition (i). Then
\begin{align*}
&\left( \frac{1}{n^{1/s}(n +1)} \right)^{s} = \frac{1}{n(n + 1)^{s}},\\[.2pc]
&I_1^s := \sum_{n = m+1}^{\infty} \frac{1}{n(n + 1)^s} \geq
\sum_{n = m+1}^{\infty}(n + 1)^{-s -1}\\[.2pc]
&\qquad\ > \int_{m+1}^{\infty}x^{-s-1}\,{\rm d}x = (1/s) (m + 1)^{-s}.
\end{align*}

Therefore condition  (ii)  of Theorem~1 takes the form of condition (ii)
of Corollary~2.

\hfill $\ab$
\end{proof}

\begin{corolla}$\left.\right.$\vspace{.5pc}

\noindent Let $\{ p_{n}\}, \{ q_{n} \}$ be positive sequences. Then 
$|\overline{N}, p_{n}|_{k} \Rightarrow | \overline{N},q_{n} |_k , k > 1$
iff\\[.5pc]

\noindent {\rm (i)}\vspace{-2.3pc}
\begin{equation*}
\frac{q_{n} P_{n}}{p_{n} Q_{n}} = O(1)
\end{equation*}
and\\[.5pc]

\noindent {\rm (ii)}\vspace{-2.85pc}
\begin{equation}
\Biggl( \sum_{n = m}^{\infty}\biggl( n^{1 - 1 / k} \frac{q_{n} }{Q_{n}
Q_{n - 1}} \biggr) ^{k} \Biggr)^{1/k} \left( \sum_{v = 1}^{m} \left|
Q_{v} - \frac{q_{v} P_{v}}{p_{v} } \right| ^{k^{\ast}}\left( \frac{1}{v}
\right) \right)^{1/k^{\ast}}\! = O(1).
\end{equation}
\end{corolla}

\begin{proof}
Corollary~3 comes from Theorem~1 by setting $s = k$. Formula (12) contains the
complicated conditions referred to on page 3 of \cite{BT92}.\hfill $\ab$
\end{proof}

\begin{corolla}$\left.\right.$\vspace{.5pc}

\noindent Let $k > 1$. Then ${\left| {C,1} \right|}_{k} \Rightarrow | 
{\overline{N},p_{n}} |_{k}$ iff\\[.5pc]

\noindent {\rm (i)}\vspace{-1.8pc}
\begin{equation}
\frac{np_{n}}{P_{n}} = O(1),
\end{equation}
\vspace{.3pc}

\noindent {\rm (ii)}\vspace{-2pc}
\begin{equation}
\frac{P_{n}}{np_{n}} = O(1)
\end{equation}
hold.
\end{corolla}

\begin{proof}
Set $s=k$. Clearly the equivalence implies (13) and (14).

To prove the converse we must show that (13) and (14) imply conditions
(ii) of Coroll- aries 1 and 2.

Using (13), with $s = k$,
\begin{equation*}
\left( {{\frac{{n^{1 - 1/k}p_{n} }}{{P_{n} P_{n - 1}} }}} \right)^{k} =
{\frac{{n^{k - 1}p_{n} ^{k}}}{{\left( {P_{n} P_{n - 1}} \right)^{k}}}} =
{\frac{{O(1)p_{n}} }{{P_{n} P_{n - 1} ^{k}}}}.
\end{equation*}

Therefore
\begin{align*}
\sum\limits_{n = m + 1}^{\infty} \left( \frac{n^{1 - 1/k }p_{n}}
{P_{n} P_{n - 1}} \right)^{k} &= O(1) \sum\limits_{n = m +
1}^{\infty} \frac{p_{n}}{P_{n} P_{n - 1}^{k}}\\[.2pc]
&= \frac{O(1)}{P_{m}^{k - 1}}\sum\limits_{n = m + 1}^{\infty}
\frac{p_{n}}{P_{n} P_{n - 1}} = \frac{O(1)}{P_{m}^{k}}.
\end{align*}

From (14),
\begin{align*}
\sum_{v = 1}^{m} \frac{|P_{v} - (v + 1)p_{v} |^{k^{\ast}}}{v} &= \sum_{v
= 1}^{m} \frac{P_{v} ^{k^{\ast}}}{v}\biggl| 1 - \frac{(v + 1)p_{v}}
{P_{v}} \biggr|^{k^{\ast}}\\[.2pc]
&= O(1)\sum_{v = 1}^{m} \frac{P_{v} ^{k^{\ast}}}{v} = O(1)\sum_{v =
1}^{m} \frac{P_{v}}{vp_{v}} p_{v} P_{v} ^{k^{\ast} - 1}
\end{align*}
\begin{align*}
&= O(1)\sum_{v = 1}^{m} p_{v} P_{v} ^{k^{\ast} - 1} \leq
O(1)P_{m}^{k^{\ast} - 1} \sum_{v = 1}^{m} p_{v}\\[.2pc]
& = O(1)P_{m} ^{k^{\ast}}.
\end{align*}

Therefore condition (ii) of Corollary 1 is satisfied.

From (14), using (13),
\begin{align*}
\sum_{v = 1}^{m} \biggl| v + 1 - \frac{P_{v}}{p_{v}}
\biggr|^{k^{\ast}}\biggl(\frac{1}{v} \biggr) & = \sum_{v =
1}^{m}\frac{(v + 1)^{k^{\ast}}}{v}\biggl| 1 - \frac{P_{v} }{(v +
1)p_{v}}\biggr|^{k^{\ast}} \\[.2pc]
&= O(1)\sum_{v = 1}^{m} v^{k^{\ast} - 1} = O(1)m^{k^{\ast}}.
\end{align*}
Therefore condition (ii) of Corollary~2 is satisfied.\hfill $\ab$
\end{proof}

Corollary~4 is Theorem~5.1 of \cite{Sa91} since (1) and (2) are
equivalent when conditions (13) and (14) hold.

We conclude by providing examples of weighted mean matrices satisfying
Corollaries~1 and 2.

\begin{exampl}
{\rm For Corollary~1, choose  $q_n = e^{-n}$. Then
$Q_n = (1 - e^{-(n+1)})/(1 - e)$, and
\begin{align*}
\frac{n^{1 + 1/k - 1/s}q_n}{Q_n} &= \frac{n^{1 + 1/k - 1/s}e^{-n}(1 -
e)}{1 - e^{-(n+1)}}\\[.2pc]
&= \frac{n^{1 + 1/k - 1/s}(1 - e)}{e^n - e^{-1}} \to 0 \quad
\textnormal{as}\quad n \to \infty,
\end{align*}
and condition (i) is satisfied.
\begin{align*}
I_2^s &:= \sum_{n = m + 1}^{\infty}\Biggl(\frac{n^{1 -
1/s}q_n}{Q_nQ_{n-1}}\Biggr)^s = \sum_{n = m + 1}^{\infty}
\Biggl(\frac{n^{1 - 1/s}e^{-n}(1 - e)^2}{(1 - e^{-(n + 1)}(1 -
e^{-n})}\Biggr)^s\\[.2pc]
&= (1 - e)^{2s}\sum_{n = m + 1}^{\infty}\frac{n^{s-1}}{(1 -
e^{-(n+1)})^s(e^n - 1)^s}\\[.2pc]
& \leq \frac{(1 - e)^{2s}}{(1 - e^{-(m+2)})^s}\sum_{n = m +
1}^{\infty}n^{s-1}e^{-ns}\biggl(\frac{e^{m+1}}{e^{m+1} - 1}\biggr)\\[.2pc]
&= O(1) \int_m^{\infty}x^{s-1}e^{-sx}\,{\rm d}x\\[.2pc]
& < O(1) \int_m^{\infty}x^{[s]}e^{-x}\,{\rm d}x\\[.2pc]
&= O(|P(m, [s])|e^{-m}),
\end{align*}
where  $P(m, [s])$ is a polynomial in  $m$  of degree $[s]$.  Therefore
\begin{align}
I_2 &= O(|P(m, [s])^{1/s}e^{-m/s}),
\end{align}
\begin{align}
I_3^{k^{\ast}} &:= \sum_{v = 1}^m|Q_v - (v +
1)q_v|^{k^{\ast}}\frac{1}{v}\nonumber\\[.2pc]
& = \sum_{v = 1}^m\frac{1}{v} \biggl|\frac{1 - e^{-(v+1)}}{1 - e} - (v +
1)e^{-v}\biggr|^{k^{\ast}}\nonumber\\[.2pc]
&= \frac{1}{(1 - e)^{k^{\ast}}} \sum_{v = 1}^m\frac{1}{v}\Bigl|1 -
e^{-(v+1)} - (v + 1)e^{-v} + (v + 1)e^{-(v+1)}\Bigr|^{k^{\ast}}\nonumber\\[.2pc]
& = O(1) \sum_{v = 1}^m\frac{1}{v} = O(\log m),\nonumber
\end{align}
and\vspace{-.6pc}
\begin{equation}
I_3 = O\biggl((\log m)^{1/k^{\ast}}\biggr).
\end{equation}

Combining (16) and (17) gives condition (ii) of Corollary~1.}
\end{exampl}

\begin{exampl}
{\rm For Corollary~2, use $p_n = 2^n$. Then $P_n = 2^{n+1} - 1$.
\begin{align*}
\frac{n^{1/k - 1/s - 1}P_n}{p_n} & = \frac{n^{1/k - 1/s - 1}(s^{n+1} -
1)}{2^n}\\[.2pc]
& =n^{1/k - 1/s - 1}(2 - 2^{-n}) \to 0 \quad \textnormal{as} \quad n \to
\infty
\end{align*}
and condition (i) is satisfied.
\begin{align*}
I_{4}^{k^{\ast}} &:= \sum_{v = 1}^m \biggl|v + 1 -
\frac{P_v}{p_v}\biggr|^{k^{\ast}}\biggl(\frac{1}{v}\biggr) = \sum_{v =
1}^m \biggl|v + 1 - \frac{(2^{v+1} -1)}{2^v}\biggr|^{k^{\ast}}\frac{1}{v}\\[.2pc]
&= \sum_{v = 1}^m \biggl|v + 1 - 2 + 2^{-v}\biggr|^{k^{\ast}}\frac{1}{v}
= \sum_{v = 1}^m \frac{1}{v}|v - 1 + 2^{-v}|^{k^{\ast}}.
\end{align*}

For $v \geq 1, 0 < v - 1 + 2^{-v} < v$. Therefore
\begin{equation*}
I_4^{k^{\ast}} < \sum_{v = 1}^mv^{k^{\ast} - 1} = O(m^{k^{\ast}}),
\end{equation*}
and condition (ii) of Corollary 2 is satisfied.}
\end{exampl}

\section*{Acknowledgement}

The first author received partial support from the Turkish Scientific
Association during the preparation of this paper.

\end{document}